\documentclass{commat}

%%% AUTHOR'S PACKAGES %%%
\usepackage{graphicx}

\title{%
    Pseudo links and singular links in the Solid Torus
    }

\author{%
    Ioannis Diamantis
    }

\affiliation{
    \address{Ioannis Diamantis --
    Department of Data Analytics and Digitalisation,
    Maastricht University, School of Business and Economics,
    P.O.Box 616, 6200 MD, Maastricht,
    The Netherlands.
    }
    \email{%
    i.diamantis@maastrichtuniversity.nl
    }
    }

\abstract{%
    In this paper we introduce and study the theories of pseudo links and singular links in the Solid Torus, ST. Pseudo links are links with some missing crossing information that naturally generalize the notion of knot diagrams, and that have potential use in molecular biology, while singular links are links that contain a finite number of self-intersections. We consider pseudo links and singular links in ST and we set up the appropriate topological theory in order to construct invariants for these types of links in ST. In particular, we formulate and prove the analogue of the Alexander theorem for pseudo links and for singular links in ST. We then introduce the mixed pseudo braid monoid and the mixed singular braid monoid, with the use of which, we formulate and prove the analogue of the Markov theorem for pseudo links and for singular links in ST.

    Moreover, we introduce the pseudo Hecke algebra of type A, $P\mathcal{H}_n$, the cyclotomic and generalized pseudo Hecke algebras of type B, $P\mathcal{H}_{1, n}$, and discuss how the pseudo braid monoid (cor. the mixed pseudo braid monoid) can be represented by $P\mathcal{H}_{n}$ (cor. by $P\mathcal{H}_{1, n}$). This is the first step toward the construction of HOMFLYPT-type invariants for pseudo links in $S^3$ and in ST. We also introduce the cyclotomic and generalized singular Hecke algebras of type B, $S\mathcal{H}_{1, n}$, and we present two sets that we conjecture they form linear bases for $S\mathcal{H}_{1, n}$. Finally, we generalize the Kauffman bracket polynomial for pseudo links in ST.
    }



\VOLUME{31}
\YEAR{2023}
\NUMBER{1}
\firstpage{333}
\DOI{https://doi.org/10.46298/cm.10438}

\begin{paper}
%\paper

\section{Introduction}\label{intro}

Pseudo diagrams of knots were introduce by Hanaki in \cite{H} as projections on the 2-sphere with over/under information at some of the double points missing. They comprise a relatively new and important model for DNA knots, since there exist DNA knots that is impossible to distinguish a positive from a negative crossing, even when studying them by electron microscopes. By considering equivalence classes of pseudo diagrams under equivalence relations generated by a specific set of Reidemeister moves, we obtain the theory of {\it pseudo knots}, that is, standard knots whose projections contain crossings with missing information (see \cite{HJMR} for more details). In \cite{BJW}, the pseudo braid monoid is introduced, which is related to the singular braid monoid, and with the use of which, the authors present the analogues of the Alexander and the Markov theorems for pseudo knots in $S^3$. From now on and throughout the paper, by pseudo links we shall mean both pseudo links and pseudo knots. In \cite{D}, the author introduces the $L$-moves for pseudo links and formulates and proves a sharpened version of the analogue of the Markov theorem for pseudo links, as well as an alternative proof of the analogue of the Alexander theorem for pseudo links in $S^3$. In this paper we extend these results for pseudo links in the solid torus ST. Our aim is the construction of pseudo knot invariants in $S^3$ and ST, and toward that end we present the pseudo bracket polynomial for pseudo links in ST and discuss HOMFLYPT-type invariants for pseudo links in $S^3$ and ST via appropriate Hecke type algebras, following \cite{Jo} \cite{La1} and \cite{PR}.

Singular knots are knots with finite many rigid self-intersections, and as shown in \cite{BJW}, the theory of singular links in $S^3$ is closely related to the theory of pseudo links in $S^3$. In this paper, we establish this correspondence for the case of ST, and through this correspondence we formulate the analogues of the Alexander and Markov theorems for singular links in ST. Then, following \cite{PR}, we define the generalized (and cyclotomic) singular Hecke algebras of type B and discuss further research needed toward the construction of HOMFLYPT-type invariants for singular links in ST, which in turn relate to invariants of pseudo links in ST. This is the subject of a sequel paper.

The paper is organized as follows: In \S~\ref{prel} we recall results on pseudo links, on singular links and on knots in the solid torus from \cite{BJW}, \cite{HJMR}, \cite{D}, \cite{La} and \cite{La2}. In \S~\ref{tpl} we present a braiding algorithm for pseudo and singular links in ST and we introduce the mixed pseudo braid monoid and the mixed singular braid monoid. We then formulate and prove the analogue of the Markov theorem for pseudo links and for singular links in ST. In \S~\ref{HOM} we introduce the pseudo Hecke algebras of type A, related to the pseudo knot theory of $S^3$, and we present some results toward the construction of HOMFLYPT-type invariants for pseudo links in $S^3$. In particular, we present a spanning set for the pseudo Hecke algebra of type A. We introduce the cyclotomic and the generalized pseudo Hecke algebras of type B, which are related to the pseudo knot theory of ST, and we discuss further steps toward the construction of HOMFLYPT-type invariants for pseudo links in ST. We conclude by presenting potential bases for these algebras. In a similar way, we also introduce the cyclotomic and generalized singular Hecke algebras of type B, through which, HOMFLYPT-type invariants for singular links in ST may be constructed. Finally, in \S~\ref{KAUF} we generalize the (Kauffman) bracket polynomial for pseudo links in ST.

The results of this paper will be used for the construction of HOMFLYPT-type invariants for these type of knots in the solid torus in the same sense as in \cite{La1} for the case of classical knots in the solid torus and in \cite{PR} for the case of singular knots in $S^3$. Our aim is to extend these invariants for the case of lens spaces $L(p,q)$ following \cite{DL3}, \cite{DL4} and \cite{DLP}.

\section{Preliminaries}\label{prel}

\subsection{The theories of pseudo links and singular links in $S^3$}\label{secpl}

In this subsection we recall results on pseudo links and on singular links in $S^3$. In particular, we recall the pseudo braid monoid $PM_n$, introduced in \cite{BJW}, that is related to the singular braid monoid $SM_n$ (\cite{Ba}, \cite{Bi}), and we recall all necessary results in order to present the analogues of the Alexander and the Markov theorems for pseudo links and singular links in $S^3$.

A {\it pseudo diagram} of a knot consists of a regular knot diagram with some missing crossing information, that is, there is no information about which strand passes over and which strand passes under the other. These undetermined crossings are called {\it pre-crossings} (for an illustration see Figure~\ref{pk1}).

\begin{figure}[ht]
\begin{center}
\includegraphics[width=1.7in]{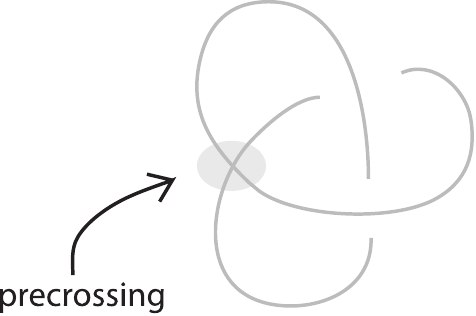}
\end{center}
\caption{A pseudo knot.}
\label{pk1}
\end{figure}

\begin{definition}
{\it Pseudo knots} are defined as equivalence classes of pseudo diagrams under an appropriate choice of Reidemeister moves that are illustrated in Figure~\ref{reid}.
\end{definition}

\begin{figure}[ht]
\begin{center}
\includegraphics[width=6.1in]{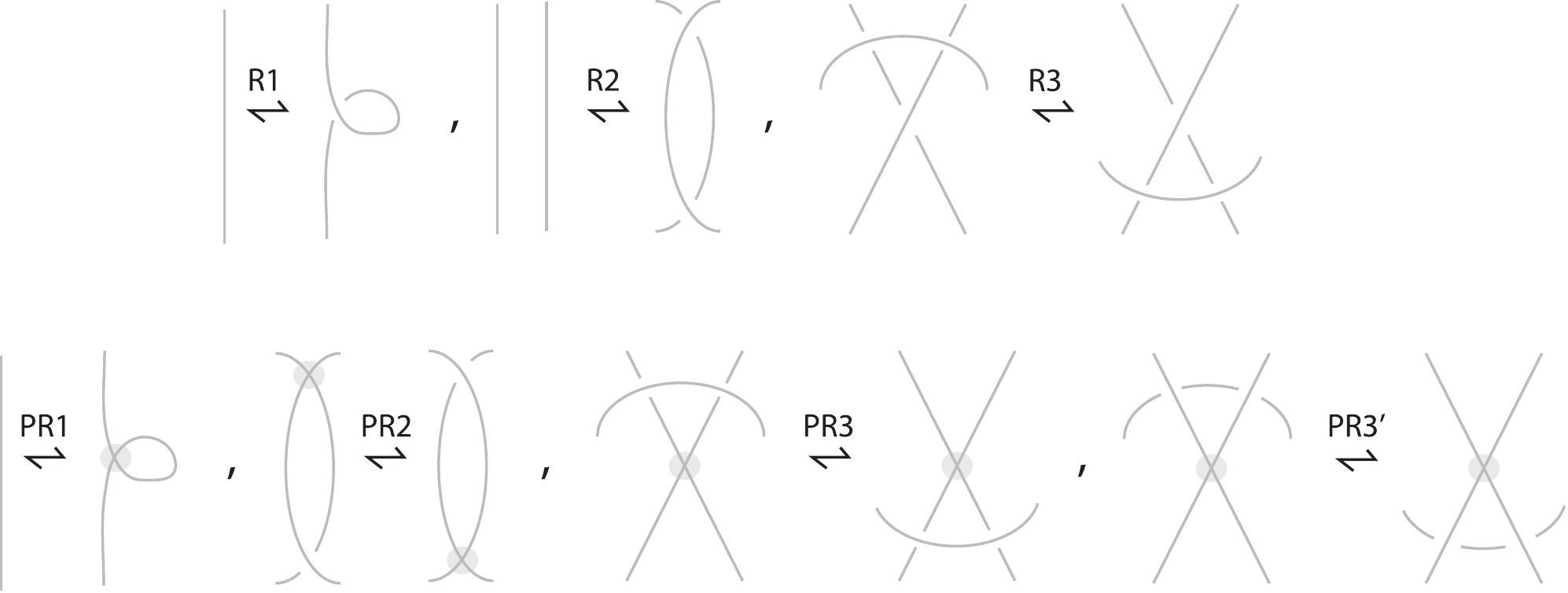}
\end{center}
\caption{Reidemeister moves for pseudo knots.}
\label{reid}
\end{figure}

As explained in \cite{BJW}, pseudo knots are closed related to {\it singular knots}, that is, knots that contain a finite number of self-intersections. In particular, there exists a bijection $f$ from the singular knot diagrams to the set of pseudo knot diagrams where singular crossings are mapped to pre-crossings. In that way we may also recover all of the pseudo knot Reidemeister moves, with the exception of the pseudo-Reidemeister I (PR1) move (see Figure~\ref{reid}). Moreover, $f$ induces an onto map from singular knots to pseudo knots, since the image of two isotopic singular knot diagrams are also isotopic pseudo knot diagrams with exactly the same sequence of Reidemeister moves.

We now introduce the pseudo braid monoid, $PM_n$, following \cite{BJW}.

\begin{definition}\label{pmn}
The monoid of pseudo braids, $PM_n$, is the monoid generated by the elements $\sigma_i^{\pm 1}, p_i$ ($i=1, \ldots, n-1$), illustrated in Figure~\ref{gens}, where $\sigma_i^{\pm 1}$ generate the braid group $B_n$ and $p_i$ satisfy the following relations:
\[
\begin{array}{rlcll}
i. & p_i\, p_j & = & p_j\, p_i, & {\rm if}\ |i-j|\geq 2\\
&&&&\\
ii. & p_i\, \sigma_j^{\pm 1} & = & \sigma_j^{\pm 1}\, p_i, & {\rm if}\ |i-j|\geq 2\\
&&&&\\
iii. & p_i\, \sigma_i^{\pm 1} & = & \sigma_i^{\pm 1}\, p_i, & i=1, \ldots, n-1\\
&&&&\\
iv. & \sigma_i\, \sigma_{i+1}\, p_i & = & p_{i+1}\, \sigma_i\, \sigma_{i+1}, & i=1, \ldots, n-2\\
&&&&\\
v. & \sigma_{i+1}\, \sigma_i\, p_{i+1} & = & p_{i}\, \sigma_{i+1}\, \sigma_i, & i=1, \ldots, n-2\\
\end{array}
\]
\end{definition}

\begin{figure}[ht]
\begin{center}
\includegraphics[width=2.7in]{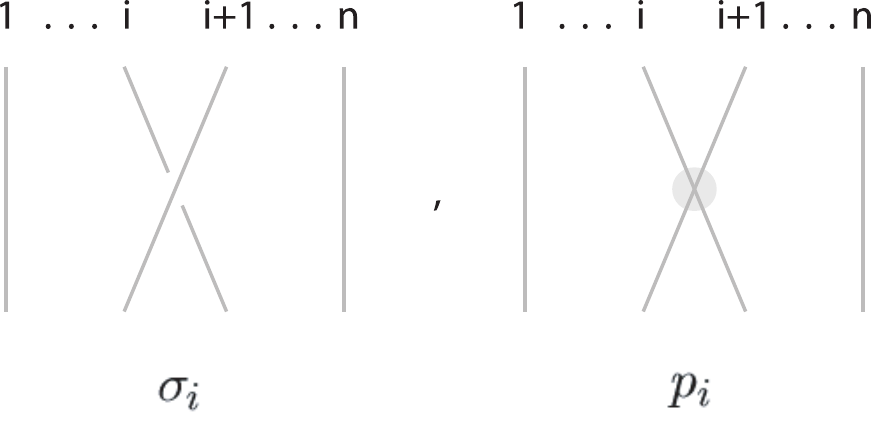}
\end{center}
\caption{The pseudo braid monoid generators.}
\label{gens}
\end{figure}

Note that $p_i$ corresponds to a standard crossing and not to a singular crossing. We denote a singular crossing by $\tau_i$ and we have that if we replace the pre-crossings $p_i$ of Definition ~\ref{pmn} by singular crossings $\tau_i$, we obtain the singular braid monoid, $SM_n$, defined in \cite{Ba} and \cite{Bi}. Thus, we obtain the following result:

\begin{proposition}[Proposition 2.3 \cite{BJW}]\label{propb1}
The monoid of pseudo braids is isomorphic to the singular braid monoid, $SM_n$. 
\end{proposition}

\begin{remark}\label{propb2}
In \cite{FKR} it is shown that $SM_n$, the singular braid monoid, embeds in a group, the singular braid group $SB_n$. It follows that $PM_n$ embeds in a group also, the pseudo braid group $PB_n$, generated by $\sigma_i$ and $p_i$, $i=1, \ldots, n-1$, satisfying the same relations as $PM_n$. Obviously, $SB_n$ is isomorphic to $PB_n$ (recall Proposition 2.3 \cite{BJW}).
\end{remark}

Define now the {\it closure} of a pseudo braid (cor. of a singular braid) as in the standard case (for an illustration of the closure of a pseudo braid see Figure~\ref{cmpl}). By considering $PM_n \subset PM_{n+1}$, we can consider the inductive limit $PM_{\infty}$. Using the analogue of the Alexander theorem for singular knots (\cite{Bi}), in \cite{BJW} the analogue of the Alexander theorem for (oriented) pseudo links is presented. In particular:

\begin{theorem}[{\bf Alexander's theorem for pseudo links}] \label{alexpl}
Every pseudo link can be obtained by closing a pseudo braid.
\end{theorem}

We now recall the notion of $L$-moves for pseudo knots and results from \cite{D} and \cite{La}. $L$-moves make up an important tool for braid equivalence in any topological setting and they allow us to formulate sharpened versions of the analogues of the Markov theorems. For more details the reader is referred to \cite{La} and references therein.

\begin{definition}\label{lmdefn}
An {\it $L$-move} on a braid $\beta$, consists in cutting an arc of $\beta$ open and pulling the upper cut-point downward and the lower upward, so as to create a new pair of braid strands with corresponding endpoints (on the vertical line of the cut-point), and such that both strands cross entirely {\it over} or {\it under} with the rest of the braid. Stretching the new strands over will give rise to an {\it $L_o$-move\/} and under to an {\it  $L_u$-move\/} as shown in Figure~\ref{lm} by ignoring all pre-crossings.
\end{definition}

$L$-moves for pseudo braids are defined in the same way as in the singular case, that is, the two strands that appear after the performance of an $L$-move should cross the rest of the braid only with real crossings (all over in the case of an $L_o$-move or all under in the case of an $L_u$-move). For an illustration see Figure~\ref{lm} for the case of an $L$-move performed on a pseudo braid. 

\begin{figure}[ht]
\begin{center}
\includegraphics[width=4.9in]{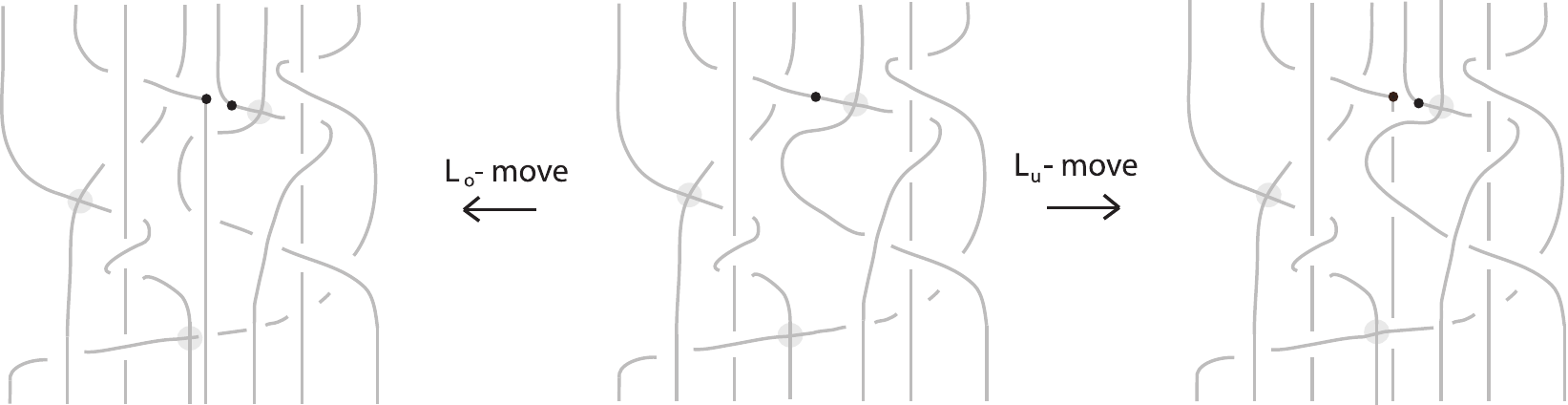}
\end{center}
\caption{$L$-moves for pseudo braids.}
\label{lm}
\end{figure}

As explained in \cite{D}, $L$-moves on pseudo knots can be used in order to obtain the analogue of the Alexander theorem for pseudo knots following the braiding algorithm in \cite{Bi}, \cite{A} and \cite{LR1}. The main idea of this algorithm is to keep the arcs of the oriented links diagram that go downwards with respect to the height function unaffected and replace arcs that go upwards with braid strands. We first isotope the pre-crossings in such a way that the braiding algorithm will not affect them (see Figure~\ref{tw}), and we may then apply the braiding algorithm of \cite{Bi} or \cite{LR1} by ignoring the pre-crossings. The same is true for the case of singular links in $S^3$ (for a detailed proof for the case of singular links in $S^3$ see \cite{PR} Theorem~2.3). In \S~\ref{secbralg} we recall the braiding algorithm from \cite{LR1} and we modify it in order to obtain the analogue of the Alexander theorem for pseudo links and singular links in ST (see also \cite{La2}).

We now state the analogue of the Markov theorem for oriented pseudo links and for oriented singular links in $S^3$ (\cite{G}, \cite{La}, \cite{BJW} and \cite{D}).

\begin{theorem}[{\bf The analogue of the Markov Theorem for singular links}]\label{mthsll}
Two oriented singular links are isotopic if and only if any two corresponding singular braids differ by braid relations in $SB_{\infty}$ and a finite sequence of the following moves:
\[
\begin{array}{lccc}
{Singular\ commuting:} & \tau_i\, \alpha & \sim & \alpha\, \tau_i\\
&&&\\
{L-moves}, & &&\\
\end{array}
\]
where $\alpha \in SB_n$, the singular braid group, and $\tau_i$ a singular crossing.
\end{theorem}

\begin{theorem}[{\bf The analogue of the Markov Theorem for pseudo links}] \label{markpl}
Two pseudo braids have isotopic closures if and only if one can be obtained from the other by a finite sequence of the following moves:
\[
\begin{array}{lllcll}
{Conjugation:} &  \alpha & \sim & \beta^{\pm 1}\, \alpha\, \beta^{\mp 1}, & {\rm for}\ \alpha \in PM_n\ \&\ \beta \in B_n,\\
&&&&\\
{Commuting:} &  \alpha\, \beta & \sim & \beta\, \alpha, & {\rm for}\ \alpha,\, \beta \in PM_n,\\
&&&&\\
{Stabilization:} &  \alpha & \sim & \alpha\, \sigma_n^{\pm 1}, & {\rm for}\ \alpha \in PM_n,\\
&&&&\\
{Pseudo-stabilization:} &  \alpha & \sim & \alpha\, p_n, & \alpha \in PM_n.\\
\end{array}
\]
Equivalently, two pseudo braids have isotopic closures if and only if one can be obtained from the other by a finite sequence of the following moves:
\[
\begin{array}{llcll}
{L-moves} &   &  &   & \\
&&&&\\
{Commuting:} &  \alpha\, \beta & \sim & \beta\, \alpha, & {\rm for}\ \alpha,\, \beta \in PM_n,\\
&&&&\\
{Pseudo-Stabilization:} &  \alpha & \sim & \alpha\, p_n, & \in PM_{n+1}.\\
\end{array}
\]
\end{theorem}

\subsection{The knot theory of ST}\label{ktst}

We now view ST as the complement of a solid torus in $S^3$ and we present results from \cite{LR1}. An oriented link $L$ in ST can be represented by an oriented \textit{mixed link} in $S^{3}$, that is, a link in $S^{3}$ consisting of the unknotted fixed part $\widehat{I}$ representing the complementary solid torus in $S^3$ and the moving part $L$ that links with $\widehat{I}$.

A \textit{mixed link diagram }is a diagram $\widehat{I}\cup \widetilde{L}$ of $\widehat{I}\cup L$ on the plane of $\widehat{I}$, where this plane is equipped with the top-to-bottom direction of $I$. For an illustration see Figure~\ref{mplink} by ignoring the pre-crossings.

Consider now an isotopy of an oriented link $L$ in ST. As the link moves in ST, its corresponding mixed link will change in $S^{3}$ by a sequence of moves that keep the oriented $\widehat{I}$ pointwise fixed. This sequence of moves consists in isotopy in $S^{3}$ and the \textit{mixed Reidemeister moves}. In terms of diagrams we have the following result for isotopy in ST:

The mixed link equivalence in $S^{3}$ includes the classical Reidemeister moves and the mixed Reidemeister moves, which involve the fixed and the standard part of the mixed link, keeping $\widehat{I}$ pointwise fixed.

By the Alexander theorem for knots in solid torus (\cite{La2}), a mixed link diagram $\widehat{I}\cup \widetilde{L}$ of $\widehat{I}\cup L$ may be turned into a \textit{mixed braid} $I\cup \beta $ with isotopic closure. This is a braid in $S^{3}$ where, without loss of generality, its first strand represents $\widehat{I}$, the fixed part, and the other strands, $\beta$, represent the moving part $L$. The subbraid $\beta$ is called the \textit{moving part} of $I\cup \beta $. For an illustration see Figure~\ref{cmpl} by ignoring the pre-crossings.

The sets of braids related to ST form groups, the Artin braid groups type B, denoted $B_{1,n}$, with presentation:

\[ B_{1,n} = \left< \begin{array}{ll}  \begin{array}{l} t, \sigma_{1}, \ldots ,\sigma_{n-1}  \\ \end{array} & \left| \begin{array}{l}
\sigma_{1}t\sigma_{1}t=t\sigma_{1}t\sigma_{1} \ \   \\
 t\sigma_{i}=\sigma_{i}t, \quad{i>1}  \\
{\sigma_i}\sigma_{i+1}{\sigma_i}=\sigma_{i+1}{\sigma_i}\sigma_{i+1}, \quad{ 1 \leq i \leq n-2}   \\
 {\sigma_i}{\sigma_j}={\sigma_j}{\sigma_i}, \quad{|i-j|>1}  \\
\end{array} \right.  \end{array} \right>, \]
where the generators $\sigma _{i}$ and $t$ are illustrated in Figure~\ref{gen}.

\begin{figure}[ht]
\begin{center}
\includegraphics[width=2.6in]{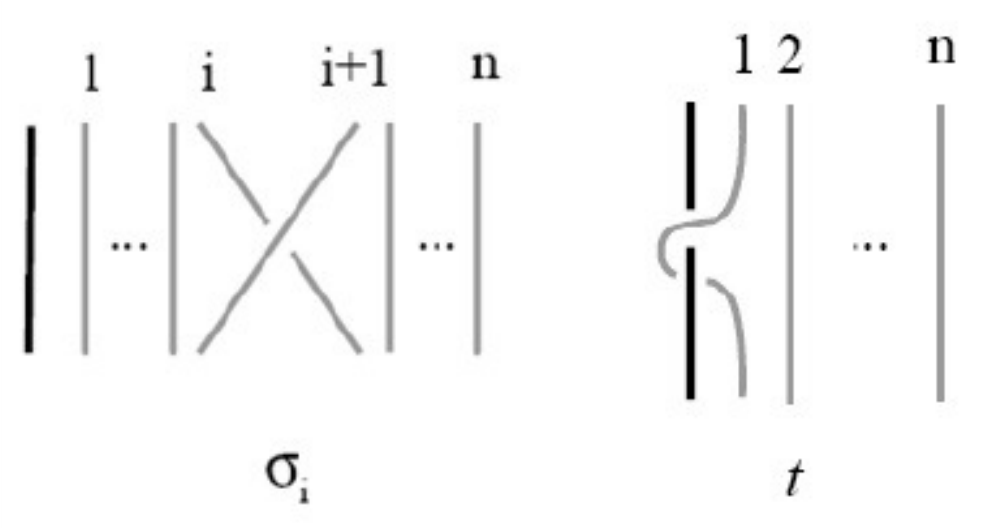}
\end{center}
\caption{The generators of $B_{1,n}$.}
\label{gen}
\end{figure}

We finally have that isotopy in ST is translated on the level of mixed braids by means of the following theorem.

\begin{theorem}[{\bf The analogue of the Markov Theorem for mixed braids}, Theorem~3 \cite{La1}] \label{markov}
 Let $L_{1} ,L_{2}$ be two oriented links in ST and let $I\cup \beta_{1} ,{\rm \; }I\cup \beta_{2}$ be two corresponding mixed braids in $S^{3}$. Then $L_{1}$ is isotopic to $L_{2}$ in ST if and only if $I\cup \beta_{1}$ is equivalent to $I\cup \beta_{2}$ in $\mathop{\cup }\limits_{n=1}^{\infty } B_{1,n}$ by the following moves:
\[ \begin{array}{clll}
(i)  & Conjugation:         & \alpha \sim \beta^{-1} \alpha \beta, & {\rm if}\ \alpha ,\beta \in B_{1,n}. \\
(ii) & Stabilization\ moves: &  \alpha \sim \alpha \sigma_{n}^{\pm 1} \in B_{1,n+1}, & {\rm if}\ \alpha \in B_{1,n}. \\
\end{array} \]
\end{theorem}

\section{Pseudo Links \& Singular Links in ST}\label{tpl}

In this section we introduce and study pseudo links and singular links in ST, through mixed pseudo (cor. singular) braids and the mixed pseudo (cor. singular) braid monoid of type B. We conclude this section by formulating and proving the analogues of the Alexander and the Markov theorems for pseudo links and singular links in ST.

\subsection{Mixed pseudo links and mixed singular links}

As explained in \S~\ref{ktst} we view ST as the complement of another ST in $S^3$ and thus, pseudo links in ST can be seen as {\it mixed pseudo links} in $S^3$ containing the complementary ST. Similarly, singular links in ST can be viewed as {\it mixed singular links} in $S^3$. For an illustration of a mixed pseudo link see Figure~\ref{mplink}, where by changing pre-crossings to singular crossings, one obtains a singular mixed link.

\begin{figure}[ht]
\begin{center}
\includegraphics[width=1.3in]{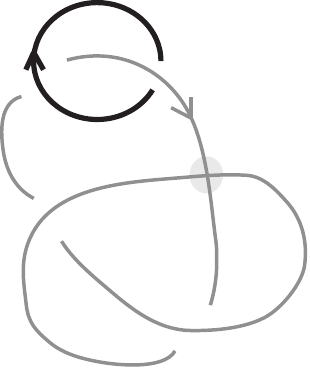}
\end{center}
\caption{A mixed pseudo link in $S^3$.}
\label{mplink}
\end{figure}

Isotopy for pseudo links in ST is now translated on the level of mixed pseudo links in $S^3$ by means of the following theorem:

\begin{theorem}[{\bf The analogue of the Reidemeister theorem for mixed pseudo links}]\label{reidplink}
Two mixed pseudo links in $S^3$ are isotopic if and only if they differ by a finite sequence of the classical and the pseudo Reidemeister moves illustrated in Figure~\ref{reid} for the standard part of the mixed pseudo links, and moves that involve the fixed and the standard part of the mixed pseudo links, called mixed Reidemeister moves, and are illustrated in Figure~\ref{mpr}.
\end{theorem}

\begin{figure}[ht]
\begin{center}
\includegraphics[width=5.1in]{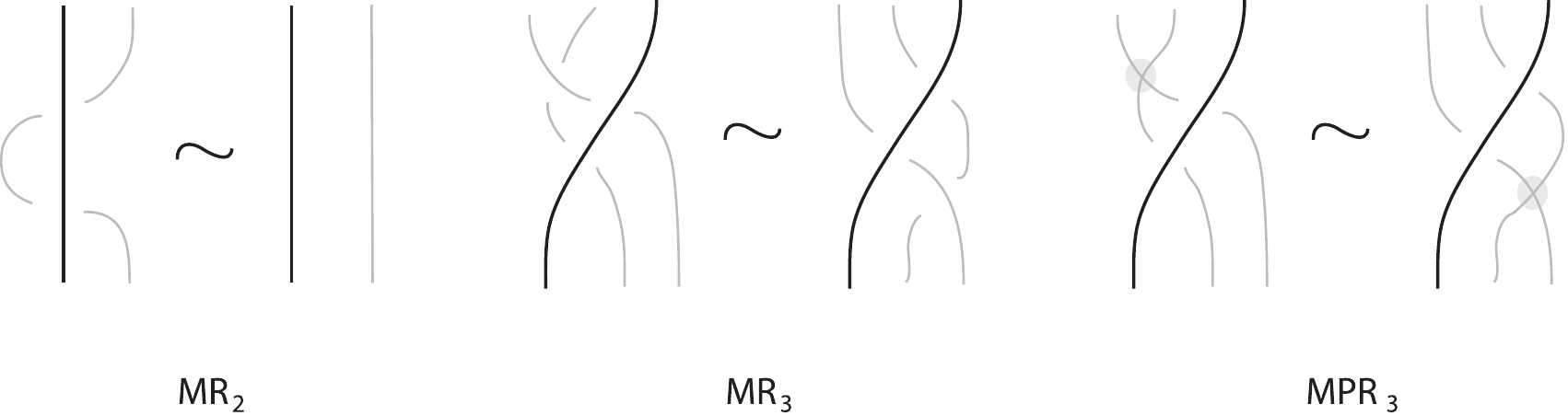}
\end{center}
\caption{Mixed Reidemeister moves.}
\label{mpr}
\end{figure}

\begin{remark}
Similarly, if we exclude PR1-moves of Figure~\ref{reid} for the standard part of the mixed singular links, and if we change the pre-crossings in Figure~\ref{mpr} to singular crossings, we obtain the analogue of the Reidemeister theorem for singular links in ST.
\end{remark}

It is worth mentioning that we do not allow {\it special} crossings, i.e. pre-crossings for the case of pseudo links and singular crossings for the case of singular links, between the fixed and the moving part of the mixed pseudo (cor. singular) braid.

\subsection{Mixed pseudo braids, mixed singular braids and a braiding algorithm}\label{secbralg}

We now define the mixed pseudo braids, which are related to mixed pseudo links, and the mixed singular braids, related to the mixed singular links.

\begin{definition}
A \textit{mixed pseudo braid} (cor. mixed singular braid) on $n$ strands, denoted by $I\cup B$, is an element of the pseudo braid monoid $PM_{1+n}$ (cor. singular braid monoid $SM_{1+n}$) consisting of two disjoint sets of strands, one of which is the identity braid $I$ representing the complementary solid torus in $S^3$, and $n$ strands form the {\it moving subbraid\/} $\beta$ representing the pseudo (cor. singular) link $L$ in ST. For an illustration see the left hand side of of Figure~\ref{cmpl}. Moreover, a diagram of a mixed pseudo braid is a braid diagram projected on the plane of $I$.
\end{definition}

Without loss of generality, we assume that the first strand of the mixed pseudo braid represents the complementary solid torus in $S^3$. This can be realized by performing the technique of (standard) parting, in order to separate the endpoints of the mixed braids into two different sets, the first represents the complementary solid torus in $S^3$ (i.e. the fixed part $I$), and the last $n$ strands represent the moving part of the mixed braid, i.e. the link in ST, and so that the resulting braids have isotopic closures. This can be realized by pulling each pair of corresponding moving strands to the right and {\it over\/} or {\it under\/} the strand of $I$ that lies on their right according to its label. We start from the rightmost pair respecting the position of the endpoints. For more details the reader is referred to \cite{LR1} and \cite{DL1}. 

We also define the closure of a mixed pseudo (cor. singular) braid as follows:

\begin{definition}
The {\it closure} $C(I \cup B)$ of a mixed pseudo braid (cor. mixed singular braid) $I \cup B$ in ST is defined as in the case of classical braids in $S^3$. See Figure~\ref{cmpl} for the case of a mixed pseudo braid.
\end{definition}

\begin{figure}[ht]
\begin{center}
\includegraphics[width=3.5in]{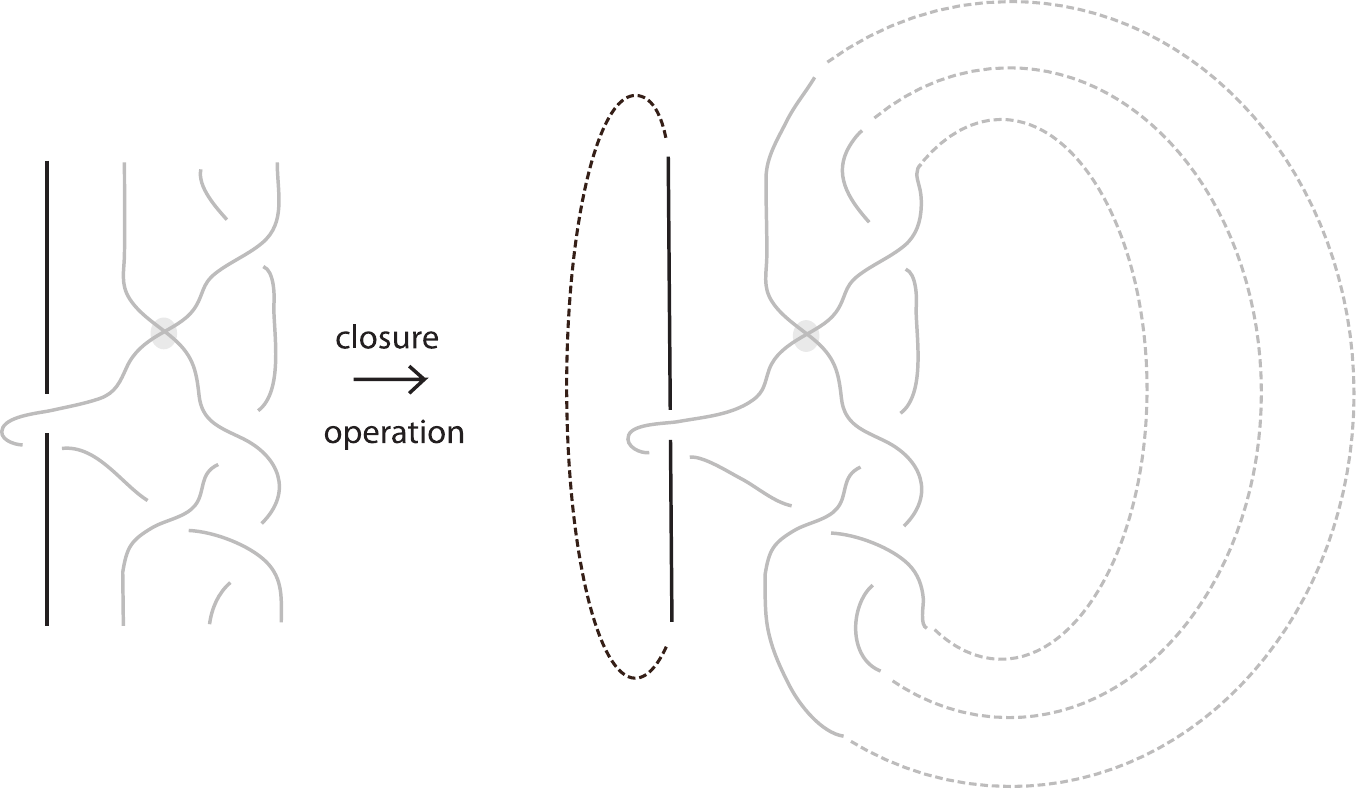}
\end{center}
\caption{The closure of a mixed pseudo braid to a mixed pseudo link.}
\label{cmpl}
\end{figure}

We now present a braiding algorithm for mixed pseudo links in ST and for mixed singular links in ST. By ``mixed links'' in the braiding algorithm that follows, we shall mean both mixed pseudo links and mixed singular links, and by {\it special crossings} we shall mean pre-crossings or singular crossings. The main idea of the braiding algorithm is to keep the arcs of the oriented mixed link diagrams that go downwards with respect to the height function unaffected, and replace arcs that go upwards with braid strands. These arcs are called {\it up-arcs} (see Figure~\ref{upa}).

\begin{figure}[ht]
\begin{center}
\includegraphics[width=4.4in]{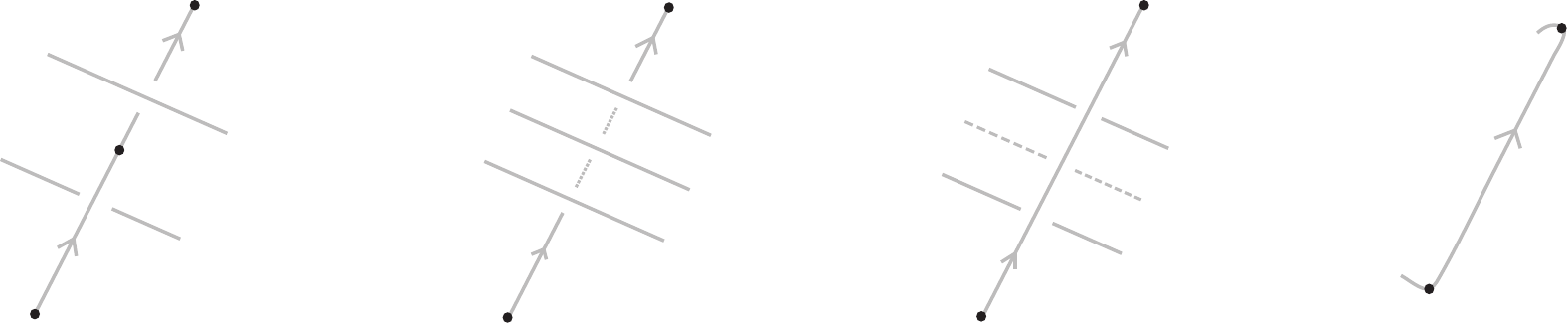}
\end{center}
\caption{Up-arcs.}
\label{upa}
\end{figure}

Note that the braiding algorithm should keep the fixed part of the mixed link unaffected and should also ''take care`` of the special crossings in the mixed diagram which contain at least one up-arc. For this we apply the idea used in \cite{Bi} for the case of singular knots, (see also \cite{KL} for the case of virtual knots). Namely, before we apply the braiding algorithm we have to isotope the mixed link in such a way that the special crossings will only contain down-arcs, so that the braiding algorithm will not affect them. This is achieved by rotating all special crossings that contain at least one up-arc, so that the two arcs are now directed downward. This is illustrated in Figure~\ref{tw} for the case of pre-crossings. Then we may apply the braiding algorithm of \cite{La2}, which is a modification of the braiding algorithm for oriented mixed links presented in \cite{LR1}, for the mixed link (ignoring the special crossings). In particular:

\begin{itemize}
\item We first ''take care`` of the special crossings by isotoping the mixed link in such a way that the braiding algorithm will not affect them. For an illustration see Figure~\ref{tw}.

\begin{figure}[ht]
\begin{center}
\includegraphics[width=5.1in]{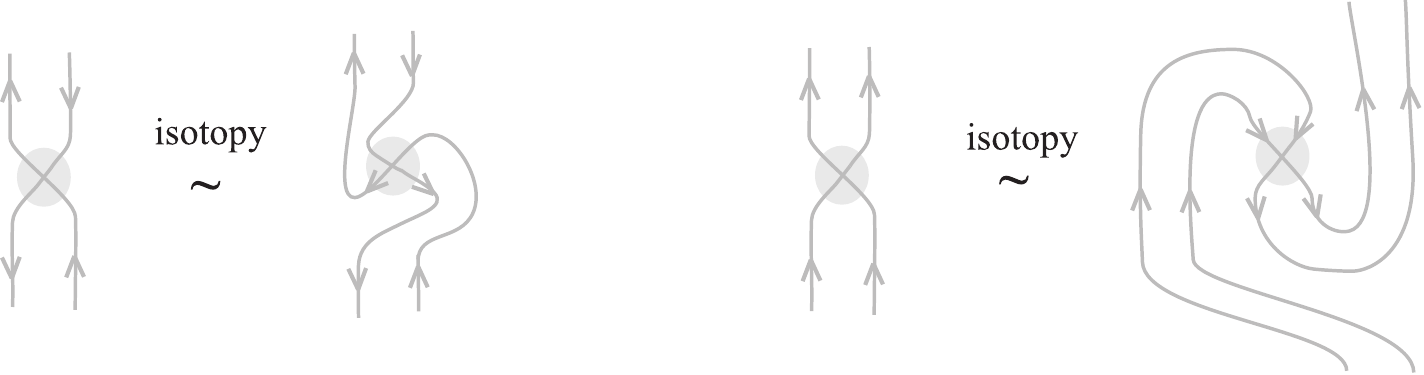}
\end{center}
\caption{Rotating pre-crossings.}
\label{tw}
\end{figure}

We now recall the braiding algorithm of \cite{La2} (for more details the reader is referred to \cite{La2} and references therein).

\item Consider now a vertical line $l$ that passes through the maximum and minimum of $\hat{I}$ as illustrated in Figure~\ref{vertl}. By small perturbations we can assume that $l$ does not pass through any crossings of the link.

\begin{figure}[ht]
\begin{center}
\includegraphics[width=2in]{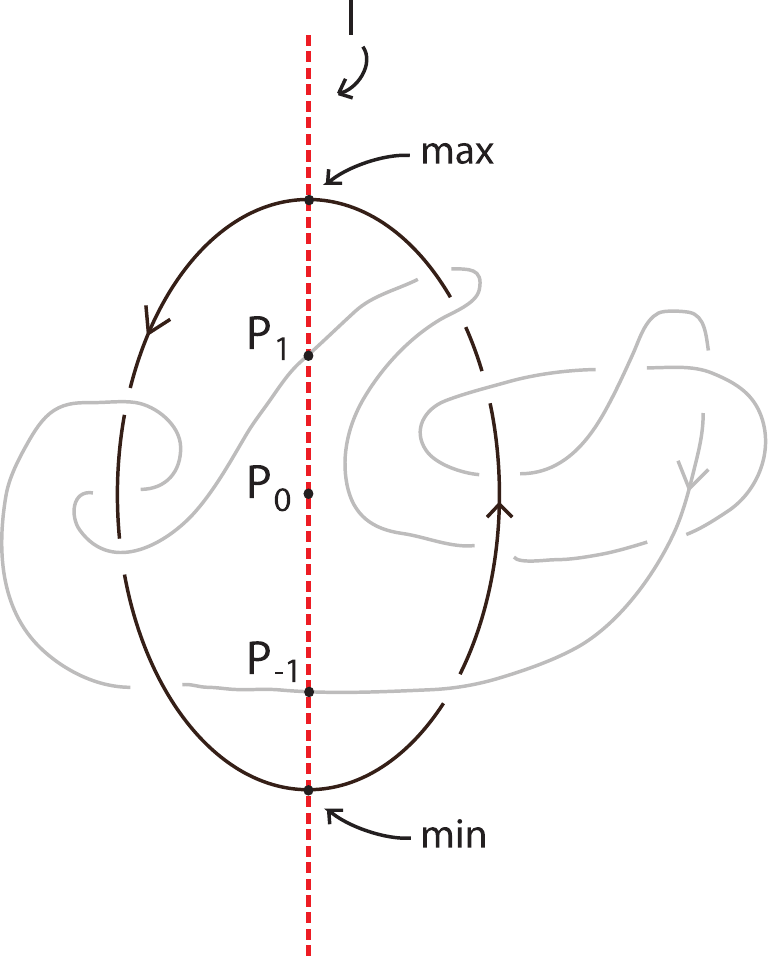}
\end{center}
\caption{}
\label{vertl}
\end{figure}

\item We now apply the following braiding algorithm to the part of the link that lies to the left of the vertical line $l$, keeping $\hat{I}$ fixed:
\begin{itemize}
\item Chose a base-point and we run along the diagram of the link according to its orientation.
 
\item When/If we run along an opposite arc, we subdivide it into smaller arcs, each containing crossings of one type only as shown in Figure~\ref{upa}.
 
\item Label now every up-arc with an ``o''or a ``u'', according to the crossings it contains. If it contains no crossings, then the choice is arbitrary.
 
\item Perform an $o$-braiding moves on all up-arcs which were labeled with an ``o'' and $u$-braiding moves on all up-arcs which were labeled with an ``u'' (see Figure~\ref{ahg}).

\begin{figure}[ht]
\begin{center}
\includegraphics[width=2.4in]{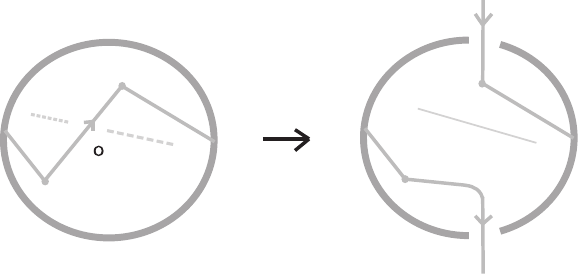}
\end{center}
\caption{Braiding moves for up-arcs.}
\label{ahg}
\end{figure}
\end{itemize}

\item Close the braided part of the link that lies to the left of $l$ and consider this operation to be enclosed in a tube $T_1$ (see left hand side of Figure~\ref{tcl}).

\item Apply once again the braiding algorithm on the right hand side of $l$ keeping $\hat{I}$ fixed, and we close now this braided portion of the link and enclose the strings participating in this operation in a tube $T_2$ (see right hand side of Figure~\ref{tcl}).

\item Rotate now around the back of the diagram in order to bring the tube $T_1$ to the very right of the diagram and $T_2$ to the very left of the diagram, and so that the resulting diagram goes around a central point, say $P_0$ on $l$.

\item Local maxima and minima in the diagram, if any, would lie on the vertical line $l$. We isolate each one in neighborhoods that do not contain other parts of the diagram, and we stretch the arcs above and below a point of symmetry $P_0$, so that the extrema would lie on $l$ again but in reverse order. Note also that we numerate the extrema with integers with respect to the point $P_0$. For an illustration see right hand side of Figure~\ref{tcl}.

\begin{figure}[ht]
\begin{center}
\includegraphics[width=4.4in]{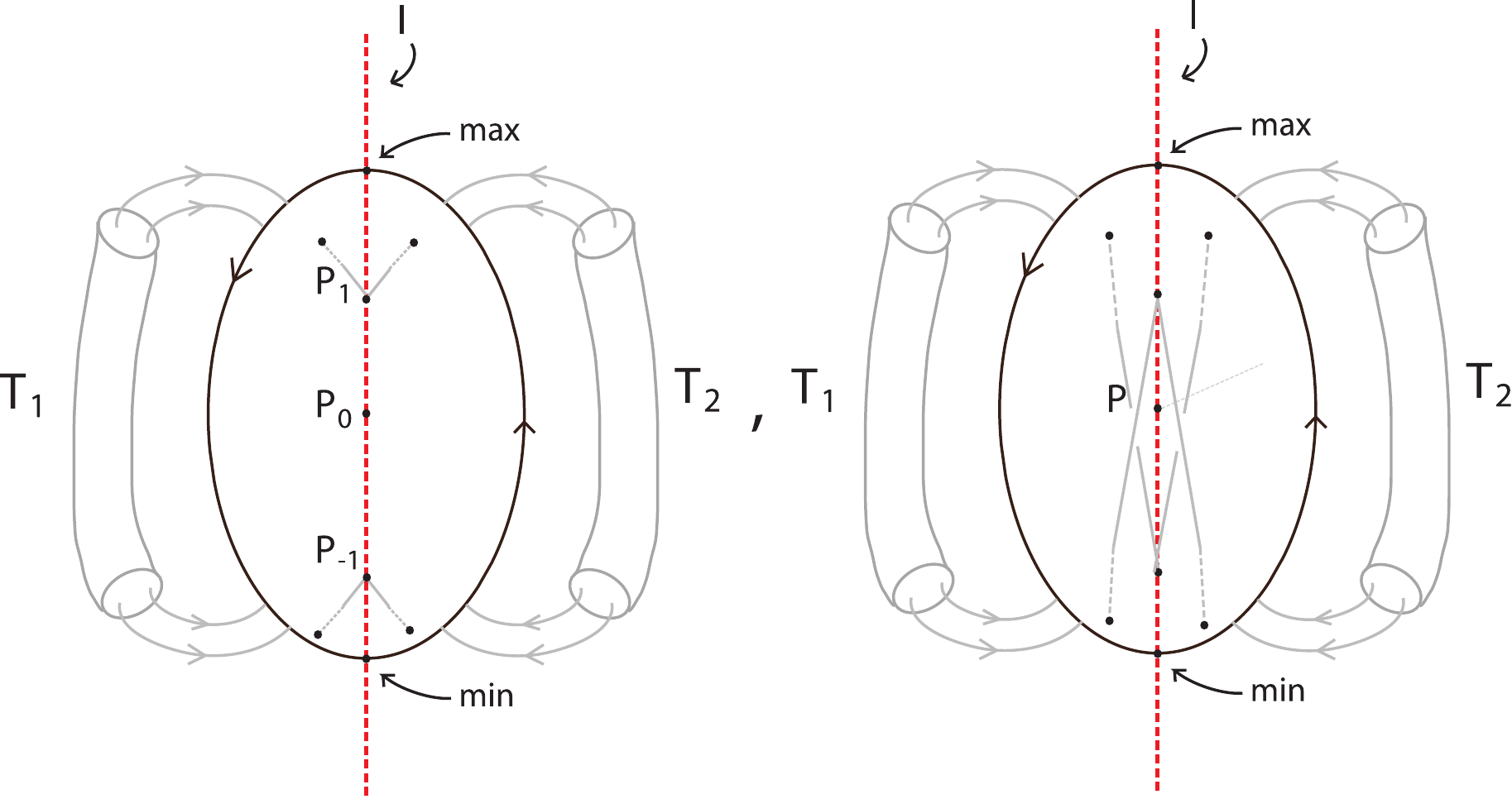}
\end{center}
\caption{}
\label{tcl}
\end{figure}

\item Open the braided diagram by cutting through a half-line starting from $P$ and we isotope in ST.

\item The result is a mixed pseudo braid (cor. a mixed singular braid) whose closure is isotopic to the initial mixed link.

\end{itemize}

This braiding algorithm provides a proof of the following theorem:

\begin{theorem}[{\bf The analogue of the Alexander theorem for oriented pseudo links and for oriented singular links in ST}]\label{newprpkalex}
Every oriented mixed pseudo link is isotopic to the closure of a mixed pseudo braid and very oriented mixed singular link is isotopic to the closure of a mixed singular braid.
\end{theorem}

\begin{remark}
\begin{itemize}
\item[i.] The braiding algorithm described above can be also applied for the case of oriented{\it tied pseudo links} in ST, that is, pseudo links in ST equipped with some non-embedded arcs called ties. The reader is referred to \cite{D} for a treatment of tied pseudo links in $S^3$. Note that tied links in $S^3$ were introduced in \cite{AJ} and that in \cite{F} the author generalized the notion of tied links in ST. Finally, in \cite{D1} the author studies tied links in various 3-manifolds.
 
\item[ii.] It is also worth mentioning that the \cite{La2}-braiding algorithm with the modifications described above, can be applied for braiding {\it pseudo singular links} in $S^3$ and ST. Pseudo singular links are defined as links with finitely many self-intersections and with some missing crossing information. This theory will be studied in a sequel paper for both $S^3$ and ST.
\end{itemize}
\end{remark}

\subsection{The mixed pseudo braid monoid \& the mixed singular braid monoid of type B}\label{secmpbm}

We now study algebraic structures related to mixed pseudo braids and mixed singular braids. Our aim is to formulate the analogue of the Markov theorem algebraically. We start by defining the {\it mixed pseudo braid monoid} $PM_{1, n}$, the counterpart of the Artin's braid group of type B for pseudo links, and also derive the {\it mixed singular braid monoid}.

In order to define the mixed pseudo braid monoid $PM_{1, n}$, we need to find the defining relations of this monoid. We analyze the moves of Theorem~\ref{reidplink} and this leads to the following definition:

\begin{definition}\label{mpbm}
The {\it mixed pseudo braid monoid of type B} $PM_{1, n}$ is defined as the monoid generated by the standard braid generators $\sigma_i^{\pm 1}$'s of $B_n$, the pseudo generators $p_{i}$'s of $PM_n$ and the looping generator $t$'s of $B_{1,n}$, satisfying the following relations:
\[
\begin{array}{rcll}
\sigma_i\, \sigma_j & = & \sigma_j\, \sigma_i, & {\rm for}\ |i-j|>1,\\
&&&\\
\sigma_i\, \sigma_j\, \sigma_i & = & \sigma_j\, \sigma_i\, \sigma_j, & {\rm for}\ |i-j|=1,\\
&&&\\
p_i\, p_j & = & p_j\, p_i, & {\rm for}\ |i-j|>1,\\
&&&\\
\sigma_i\, p_j & = & p_j\, \sigma_i, & {\rm for}\ |i-j|>1,\\
&&&\\
\sigma_i\, \sigma_j\, p_i & = & p_j\, \sigma_i\, \sigma_j, & {\rm for}\ |i-j|=1,\\
&&&\\
t\, \sigma_i & = & \sigma_i\, t, & {\rm for}\ i>1,\\
&&&\\
t\, p_i & = & p_i\, t, & {\rm for}\ i>1,\\
&&&\\
t\, \sigma_1\, t\, \sigma_1 & = & \sigma_1\, t\, \sigma_1\, t, &\\
&&&\\
t\, \sigma_1\, t\, p_1 & = & p_1\, t\, \sigma_1\, t, &\\
&&&\\
\sigma_1\, \sigma_1^{-1} & = & t\, t^{-1}\ =\ 1.&
\end{array}
\]
\end{definition}

Note that every isotopy can be decomposed in a sequence of elementary isotopies which correspond to the relations in Definition~\ref{mpbm}, and that means that we have obtained a complete set of relations. It is also worth mentioning that by replacing the pre-crossings $p_i$ in Definition~\ref{mpbm}, we obtain the {\it singular braid monoid of type B}, $SM_{1, n}$, presented in Theorem~4.4 \cite{V}. Note also that in \cite{V}, the singular braid monoid of type B is called singular braid monoid of the Annulus. Thus, we have the following result:

\begin{theorem}\label{is1}
There exists an isomorphism $\mu$ from the singular braid monoid of type B to the pseudo braid monoid of type B, defined as follows:
\begin{equation}\label{iso}
\begin{array}{rccl}
\mu\, : & SM_{1, n} & \rightarrow & PM_{1, n}\\
&&&\\
        &  \sigma_i^{\pm 1} & \mapsto & \sigma_i^{\pm 1}\\
&&&\\
	      &  t^{\pm 1}        & \mapsto & t^{\pm 1}\\
&&&\\
				&  \tau_i           & \mapsto & p_i\\
\end{array}
\end{equation}
\end{theorem}

Moreover, by considering the natural inclusion $PM_{1, n} \subset PM_{1, n+1}$ (see Figure~\ref{dirlim}), we can define the inductive limit $PM_{1, \infty}$, and similarly, by considering the natural inclusion $SM_{1, n} \subset SM_{1, n+1}$, we may define the inductive limit $SM_{1, \infty}$.

\begin{figure}[ht]
\begin{center}
\includegraphics[width=3.4in]{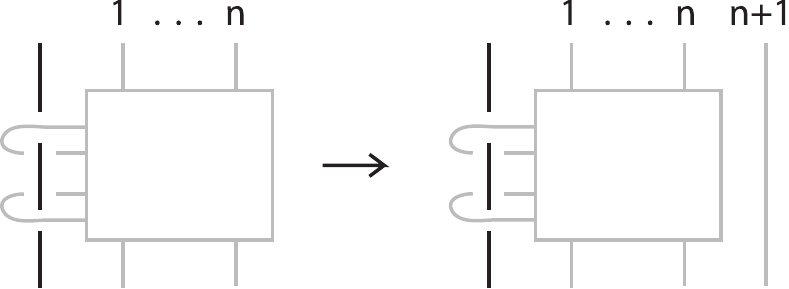}
\end{center}
\caption{The natural inclusion $PM_{1, n} \subset PM_{1, n+1}$.}
\label{dirlim}
\end{figure}

\subsection{The analogue of the Markov theorem for mixed pseudo braids and mixed singular braids}\label{mthmtpb}

Using the braiding algorithm presented in \S~\ref{secbralg} for pseudo links in $S^3$, we have a stronger version of the analogue of the Markov theorem for pseudo knots, similar to Theorem~2 in \cite{La2}. In particular, we have the following:

\begin{theorem}[{\bf Relative version of the analogue of the Markov theorem for pseudo links}]\label{rvmpl}
Two pseudo links containing the same braided part are isotopic if and only if any two corresponding pseudo braids, both containing the same braided part, differ by conjugation, commuting, stabilization and pseudo-stabilization moves that do not affect the already braided part.
\end{theorem}

Since all pseudo links related to the solid torus contain the same braided part $\hat{I}$, we have that corresponding pseudo braids contain the same braided part $I$. Therefore, we obtain the following:

\begin{corollary}\label{crvmpl}
Two pseudo links in ST are isotopic if and only if any two corresponding pseudo braids of theirs, differ by conjugation, commuting, stabilization and pseudo-stabilization moves that do not affect $I$.
\end{corollary}

By Corollary~\ref{crvmpl} and the discussion in \S~\ref{secmpbm} we obtain the following result:

\begin{theorem}[{\bf The analogue of the Markov Theorem for mixed pseudo braids}] \label{marktpb}
Two mixed pseudo braids have equivalent closures if and only if one can  obtained from the other by a finite sequence of the following moves:
\[
\begin{array}{llcll}
{\rm Commuting:} &  \alpha\, p_i & \sim & p_i\, \alpha, & {\rm for\ all}\ \alpha \in PM_{1, n},\\
&&&&\\
{\rm Conjugation:} & \beta & \sim & \alpha^{\pm 1}\, \beta\, \alpha^{\mp 1} & {\rm for\ all}\ \beta\in PM_{1, n}\ \&\ \alpha\in B_{1, n},\\
&&&&\\
{\rm Real-Stabilization:} &  \alpha & \sim & \alpha\, \sigma_n^{\pm 1}, & {\rm for\ all}\ \alpha \in PM_{1, n},\\
&&&&\\
{\rm Pseudo-Stabilization:} &  \alpha & \sim & \alpha\, p_n, & {\rm for\ all}\ \alpha \in PM_{1, n}.\\
\end{array}
\]
\end{theorem}

Note that from the discussion in \S~\ref{secpl} and Theorem~\ref{markpl} it follows that conjugation and stabilization moves can be replaced by the $L$-moves. Hence, we obtained a sharpened version of Theorem~\ref{marktpb}. More precisely, we have the following:

\begin{theorem}[{\bf $L$-move Markov's Theorem for mixed pseudo braids}] \label{lmarkmpl}
Two mixed pseudo braids have isotopic closures if and only if one can  obtained from the other by a finite sequence of the following moves:
\[
\begin{array}{llcll}
{\rm L-moves} &   &  &   & \\
&&&&\\
{\rm Commuting:} &  \alpha\, p_i & \sim & p_i\, \alpha, & {\rm for\ all}\ \alpha \in PM_{1, n},\\
&&&&\\
{\rm Pseudo-Stabilization:} &  \alpha & \sim & \alpha\, p_n, & {\rm for\ all}\ \alpha \in PM_{1, n}.\\
\end{array}
\]
\end{theorem}

Similarly, by Theorem~\ref{is1}, we have the following result for mixed singular braids:

\begin{theorem}[{\bf The analogue of the Markov Theorem for mixed singular braids}] \label{marktsb}
Two mixed singular braids have equivalent closures if and only if one can  obtained from the other by a finite sequence of the following moves:
\[
\begin{array}{llcll}
{\rm Commuting:} &  \alpha\, p_i & \sim & p_i\, \alpha, & {\rm for\ all}\ \alpha \in SM_{1, n},\\
&&&&\\
{\rm Conjugation:} & \beta & \sim & \alpha^{\pm 1}\, \beta\, \alpha^{\mp 1} & {\rm for\ all}\ \beta\in SM_{1, n}\ \&\ \alpha\in B_{1, n},\\
&&&&\\
{\rm Real-Stabilization:} &  \alpha & \sim & \alpha\, \sigma_n^{\pm 1}, & {\rm for\ all}\ \alpha \in SM_{1, n},\\
\end{array}
\]
Equivalently, by a finite sequence of the following moves:
\[
\begin{array}{llcll}
{\rm L-moves} &   &  &   & \\
&&&&\\
{\rm Commuting:} &  \alpha\, p_i & \sim & p_i\, \alpha, & {\rm for\ all}\ \alpha \in SM_{1, n},\\
\end{array}
\]
\end{theorem}

\section{Pseudo Hecke and singular Hecke algebras of type A and of type B}\label{HOM}

In \cite{PR} the authors define the singular Hecke algebra of type A, $S\mathcal{H}_n$, and through a universal Markov trace constructed in $S\mathcal{H}_n$, they present a HOMFLYPT-type invariant for singular links in $S^3$. In this section we discuss techniques from \cite{Jo}, \cite{PR} and \cite{La1} toward the construction of HOMFLYPT-type invariants for pseudo links in $S^3$ and ST, and for singular links in ST. The construction of these type of invariants will be the subject of a sequel paper.

\subsection{The pseudo Hecke algebra of type A}

The technique of \cite{PR} can be applied for the case of pseudo knots in $S^3$, since, as mentioned before, pseudo knot theory of $S^3$ can be realized as the quotient of the theory of singular knots in $S^3$, modulo the pseudo-Reidemeister move 1. In that sense, we present here results from \cite{PR} which are adapted accordingly.

Recall that the Hecke algebra of type A, $\mathcal{H}_n$, has presentation obtained from the classical braid group $B_n$ by corresponding the braiding generators $\sigma_i$ to $g_i$ and by adding the quadratic relations:
\begin{equation}\label{quadr}
g_i^2\ =\ (q-1)\, g_i\ +\ q, \qquad 1\leq i < n,\ q\in \mathbb{C}.
\end{equation}
That is, $\mathcal{H}_n\ =\ \frac{\mathbb{C}[B_n]}{\langle g_i^2-(q-1)g_i-q \rangle}$.

We define the pseudo Hecke algebra of type A in a similar way. More precisely:

\begin{definition}
Define the {\it pseudo Hecke algebra of type A}, $P\mathcal{H}_n$, as the quotient of the pseudo braid monoid $\mathbb{C}[PM_n]$ by the quadratic relations (\ref{quadr}). That is,
\[
P\mathcal{H}_n\ =\ \frac{\mathbb{C}[PM_n]}{\langle g_i^2-(q-1)g_i-q \rangle}.
\]
\end{definition}

Note that the embedding $PM_n \hookrightarrow PM_{n+1}$ induces a homomorphism $\iota : P\mathcal{H}_n \to P\mathcal{H}_{n+1}$. Moreover, each $PM_n$ has a natural grading with respect to the pre-crossings, namely, $PM_n = \bigoplus_{d=0}^\infty P_dM_n$, where $P_dM_n$ denote the set of mixed pseudo braids with $d$ pre-crossings. This grading induces a grading on $P\mathcal{H}_n$, since the quadratic relations affect only the $g_i$'s. In particular, we have that:

\[
P\mathcal{H}_n\ =\ \underset{d=0}{\overset{\infty}{\oplus}}\, P_d\mathcal{H}_n.
\]

Recall now the following linear basis for the Hecke algebra of type A \cite{Jo}:

$$S =\left\{(g_{i_{1} }g_{i_{1}-1}\ldots g_{i_{1}-k_{1}})(g_{i_{2} }g_{i_{2}-1 }\ldots g_{i_{2}-k_{2}})\ldots (g_{i_{p} }g_{i_{p}-1 }\ldots g_{i_{p}-k_{p}})\right\},$$
for $1\le i_{1}<\ldots <i_{p} \le n-1{\rm \; }$.

The basis $S$ yields directly an inductive basis for $\mathcal{H}_{n}$, which is used in the construction of the Ocneanu trace, leading to the HOMFLYPT or $2$-variable Jones polynomial. Similarly, we want to find a basis for $P\mathcal{H}_n$. Although $P\mathcal{H}_n$ is not of finite dimension, each subspace $P_d\mathcal{H}_n$ of the graduation is of finite dimension. Indeed we have the following (compare to Proposition~3.1 \cite{PR}):

\begin{proposition}\label{bas1}
The set 
\[
C_d\ =\ \left\{ p_{i_1}\, p_{i_2}\, \ldots\, p_{i_d}\, \alpha,\ {\rm where}\ 1\leq i_j\leq n-1,\, 1\leq j \leq d\ {\rm and}\ \alpha\in B_n\right\},
\]
spans $P_d\mathcal{H}_n$.
\end{proposition}

The proof is similar to that of Proposition~3.1 \cite{PR} by replacing the singular crossings $\tau_i$'s by the pre-crossings $p_i$'s.

In order now to define HOMFLYPT-type invariants for pseudo links in $S^3$, one should construct a family of Markov traces on $\{P_d\mathcal{H}_n\}_{n=1}^{\infty}$ and $P\mathcal{H}_n$. Note that the difference between the family of traces we have to define in $\{P_d\mathcal{H}_n\}_{n=1}^{\infty}$ to that defined in \cite{PR}, is the relation that comes from the pseudo-Reidemeister move 1, and which is not allowed in the theory of singular knots.

In a sequel paper we shall work toward the construction of such Markov traces. Note that the construction of a Markov trace leads naturally to a HOMFLYPT-type invariant for pseudo links in $S^3$. As shown in \cite{PR}, one would have to fix a Markov trace on $\{P_d\mathcal{H}\}_{n=1}^{\infty}$, and then follow \cite{Jo} in order to define an invariant $X$ for pseudo knots in $S^3$.

\subsection{On the pseudo Hecke algebras and the singular Hecke algebras of type B}

In this subsection we define the {\it generalized pseudo Hecke algebra of type B} and the {\it cyclotomic pseudo Hecke algebras of type B}, related to pseudo links in ST, and similarly, we define the {\it cyclotomic and the generalized singular Hecke algebra of type B}, related to the theory of singular links in ST, following \cite{La1}. Our aim is to apply similar techniques from \cite{La1}, \cite{PR} and \cite{Jo} toward the construction of HOMFLYPT-type invariants for pseudo links and singular links in the Solid Torus, ST, which is the subject of a sequel paper.

It has been established that Hecke algebras of type B form a tower of B-type algebras that are related to the knot theory of ST. A presentation for the basic one is obtained from the presentation of the Artin group $B_{1,n}$ by adding the quadratic relation~(\ref{quadr}) and $t^{2} =\left(Q-1\right)t+Q$, where $q,Q \in {\mathbb C}\backslash \{0\}$ are seen as fixed variables. The middle B-type algebras are the cyclotomic Hecke algebras of type B, whose presentations are obtained by the quadratic relation~(\ref{quadr}) and $t^d=(t-u_{1})(t-u_{2}) \ldots (t-u_{d})$. The topmost Hecke-like algebra in the tower is the \textit{generalized Hecke algebra of type B}, $\textrm{H}_{1,n}$, which has the following presentation:
\[
\textrm{H}_{1,n}{(q)} = \left< \begin{array}{ll}  \begin{array}{l} t, g_{1}, \ldots ,g_{n-1}  \\ \end{array} & \left| \begin{array}{l} g_{1}tg_{1}t=tg_{1}tg_{1} \ \
\\
 tg_{i}=g_{i}t, \quad{i>1}  \\
{g_i}g_{i+1}{g_i}=g_{i+1}{g_i}g_{i+1}, \quad{1 \leq i \leq n-2}   \\
 {g_i}{g_j}={g_j}{g_i}, \quad{|i-j|>1}  \\
 {g_i}^2=(q-1)g_{i}+q, \quad{i=1,\ldots,n-1}
\end{array} \right.  \end{array} \right>.
\]
That is:
\begin{equation*}
\textrm{H}_{1,n}(q)= \frac{{\mathbb C}\left[q^{\pm 1} \right]B_{1,n}}{ \langle \sigma_i^2 -\left(q-1\right)\sigma_i-q \rangle}.
\end{equation*}

Note that in $\textrm{H}_{1,n}(q)$ the generator $t$ satisfies no polynomial relation, making the algebra $\textrm{H}_{1,n}(q)$ infinite dimensional. Also that in \cite{La1} the algebra $\textrm{H}_{1,n}(q)$ is denoted as $\textrm{H}_{n}(q, \infty)$.

Following \cite{La1} and \cite{PR}, we now define the generalized pseudo Hecke algebra of type B and the cyclotomic pseudo Hecke algebras of type B.

\begin{definition}
Define the {\it generalized pseudo Hecke algebra of type B}, $P\mathcal{H}_{1, n}$, as the quotient of the mixed pseudo braid monoid $\mathbb{C}[PM_{1, n}]$ by the quadratic relations (\ref{quadr}). That is,
\[
P\mathcal{H}_{1, n}\ =\ \frac{\mathbb{C}[PM_{1, n}]}{\langle g_i^2-(q-1)g_i-q \rangle}.
\]

The cyclotomic pseudo Hecke algebras of type B, $P\mathcal{H}^c_{1, n}$ is defined as the quotient of the mixed pseudo braid monoid $\mathbb{C}[PM_{1, n}]$ by the quadratic relations (\ref{quadr}) and the relations:
\[
t^c=(t-u_{c})(t-u_{c}) \ldots (t-u_{c}).
\]
\end{definition}

As in the case of pseudo Hecke algebras of type A, the embedding $PM_{1, n} \hookrightarrow PM_{1, n+1}$ induces a homomorphism $\iota:P\mathcal{H}_{1, n}\rightarrow P\mathcal{H}_{1, n+1}$, and moreover, $PM_{1, n}$ also has a natural grading with respect to the pre-crossings, i.e. $PM_{1, n} = \bigoplus_{d=0}^\infty P_dM_{1, n}$, where $P_dM_{1, n}$ denotes the set of mixed pseudo braids with $d$ pre-crossings. This grading induces a grading on $P\mathcal{H}_{1, n}$:

\[
P\mathcal{H}_{1, n}\ =\ \underset{d=0}{\overset{\infty}{\oplus}}\, P_d\mathcal{H}_{1, n}.
\]

Consider now the elements illustrated in Figure~\ref{genh} and recall the following linear bases for the generalized Hecke algebra of type B \cite{La1}:

\[
\begin{array}{llll}
 (i) & \Sigma_{n} & = & \{t_{i_{1} } ^{k_{1} } \ldots t_{i_{r}}^{k_{r} } \cdot \sigma \} ,\ {\rm where}\ 0\le i_{1} <\ldots <i_{r} \le n-1,\\
 (ii) & \Sigma^{\prime} _{n} & = & \{ {t^{\prime}_{i_1}}^{k_{1}} \ldots {t^{\prime}_{i_r}}^{k_{r}} \cdot \sigma \} ,\ {\rm where}\ 0\le i_{1} < \ldots <i_{r} \le n-1, \\
\end{array}
\]
where $k_{1}, \ldots ,k_{r} \in {\mathbb Z}$ and $\sigma$ a basic element in ${\rm H}_{n}(q)$.

\begin{figure}[ht]
\includegraphics[width=3.5in]{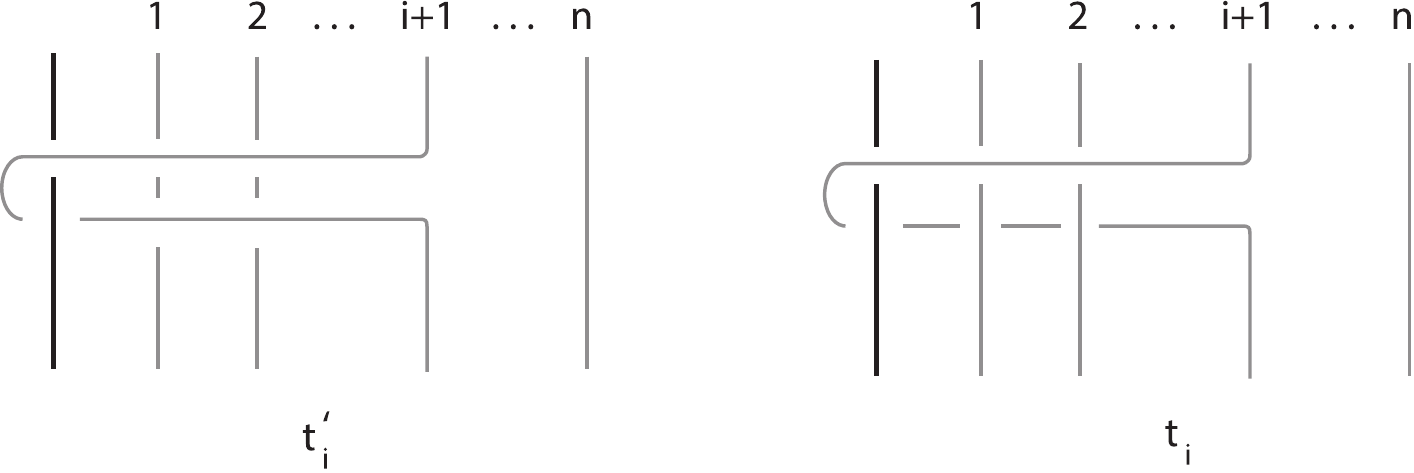}
\caption{The elements $t^{\prime}_{i}$ and $t_{i}$.}
\label{genh}
\end{figure}

The basis $\Sigma^{\prime} _{n}$ yields an inductive basis for the generalized Hecke algebra of type B, which is used in the construction of the Ocneanu trace, leading to the HOMFLYPT or $2$-variable Jones polynomial (the reader is referred to \cite{La1} for more details). Similarly, we want to find a basis for $P\mathcal{H}_{1, n}$, with the use of which we will construct HOMFLYPT-type invariants for pseudo links in ST. This is more complicated (compared to the type A case) and further research is required.

We conclude this discussion with a conjecture:

\begin{conjecture}\label{c1}
For $\sigma\in P_d\mathcal{H}_n$, the following sets span $P_d\mathcal{H}_{1, n}$:
\[
\begin{array}{llll}
 (i) & C^d_{n} & = & \{t_{i_{1} } ^{k_{1} } \ldots t_{i_{r}}^{k_{r} } \cdot \sigma \} ,\ {\rm where}\ 0\le i_{1} <\ldots <i_{r} \le n-1,\\
 (ii) & {C^{\prime}}^d _{n} & = & \{ {t^{\prime}_{i_1}}^{k_{1}} \ldots {t^{\prime}_{i_r}}^{k_{r}} \cdot \sigma \} ,\ {\rm where}\ 0\le i_{1} < \ldots <i_{r} \le n-1. \\
\end{array}
\]
\end{conjecture}

\begin{remark}
For singular links in ST, we define the {\it generalized singular Hecke algebra of type B}, $S\mathcal{H}_{1, n}$, as

\[
S\mathcal{H}_{1, n}\ =\ \frac{\mathbb{C}[SM_{1, n}]}{\langle g_i^2-(q-1)g_i-q \rangle},
\]
and the {\it cyclotomic singular Hecke algebras of type B}, $S\mathcal{H}^c_{1, n}$, as the following quotient:

\[
S\mathcal{H}^c_{1, n}\ =\ \frac{\mathbb{C}[SM_{1, n}]}{\langle g_i^2-(q-1)g_i-q,\ t^c=(t-u_{c})(t-u_{c}) \ldots (t-u_{c}) \rangle}.
\]

Similarly to Conjecture~\ref{c1}, the sets
\[
\begin{array}{llll}
 (i) & C^d_{n} & = & \{t_{i_{1} } ^{k_{1} } \ldots t_{i_{r}}^{k_{r} } \cdot \sigma \} ,\ {\rm where}\ 0\le i_{1} <\ldots <i_{r} \le n-1,\\
 (ii) & {C^{\prime}}^d _{n} & = & \{ {t^{\prime}_{i_1}}^{k_{1}} \ldots {t^{\prime}_{i_r}}^{k_{r}} \cdot \sigma \} ,\ {\rm where}\ 0\le i_{1} < \ldots <i_{r} \le n-1. \\
\end{array},
\]
where $\sigma\in S_d\mathcal{H}_n$, are potential spanning sets of $S_d\mathcal{H}_{1, n}$.
\end{remark}

\begin{remark}
It is worth mentioning that, as shown in \cite{DL3}, \cite{DL4}, \cite{DLP} and \cite{D4}, the generalized Hecke algebra of type B is related to the knot theory of lens spaces $L(p, q)$. This can be realized by the fact that one may generalize a HOMFLYPT-type invariant for links in ST to an invariant for  links in $L(p, q)$ by imposing relations coming from the {\it braid band moves}, that is, moves that reflect isotopy in $L(p, q)$ and which are similar to the second Kirby move (\cite{LR1}, \cite{DL1}).
\end{remark}

\section{The pseudo bracket polynomial for pseudo links in ST}\label{KAUF}

In this section we recall results from \cite{LK} and \cite{HD} on the pseudo bracket polynomial for pseudo links in $S^3$ and we extend this polynomial for the case of pseudo links in ST. For that reason, we view ST as a punctured disk, that is, a disk with a hole in its center, representing the complementary solid torus in $S^3$ (for an illustration see Figure~\ref{pstt}).

\subsection{The pseudo bracket polynomial for pseudo links in $S^3$}

In \cite{HD} the pseudo bracket polynomial, $\langle ; \rangle$, is defined for pseudo links in $S^3$ extending the Kauffman bracket polynomial for classical knots presented in \cite{LK}. Note that the orientation of a diagram in the case of pseudo links is needed in order to define a skein relation on pre-crossings. More precisely, we have the following:

\begin{definition}\label{pkaufb}
Let $L$ be an oriented pseudo link in $S^3$. The {\it pseudo bracket polynomial} of $L$ is defined by means of the following relations:
\begin{figure}[ht]
\begin{center}
\includegraphics[width=3.6in]{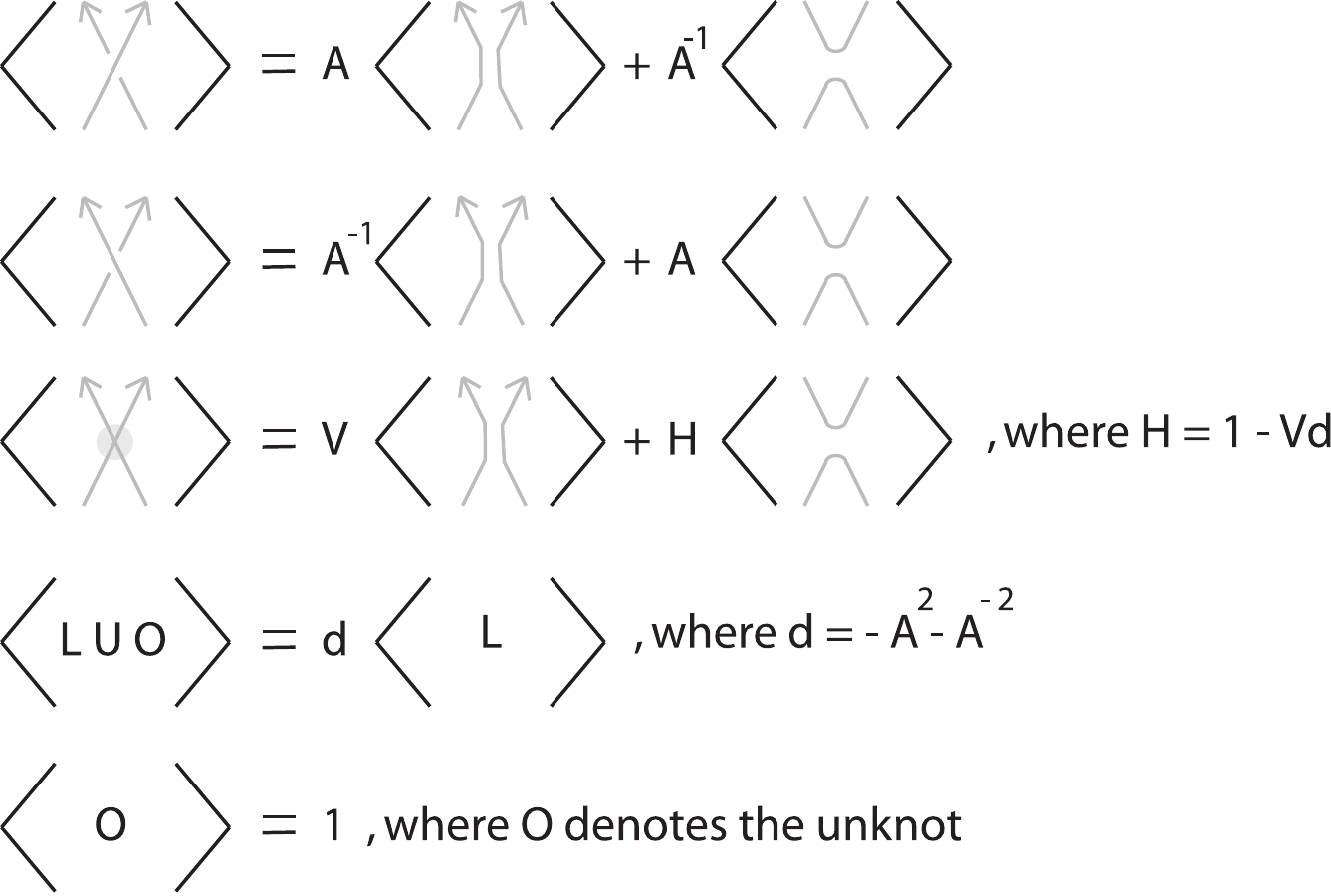}
\end{center}
%\caption{}
\label{pkb}
\end{figure}
\end{definition}

It can be easily seen that the pseudo bracket is invariant under Reidemeister moves 2 and 3 and the pseudo moves PR1, PR2 and PR3 (\cite{HD} Theorem~1). As in \cite{LK}, by normalizing the pseudo bracket polynomial using the writhe, we obtain the {\it normalized pseudo bracket polynomial}, which is an invariant of pseudo knots in $S^3$ (\cite{HD}, Corollary~2). In particular:

\begin{theorem}
Let $K$ be a pseudo diagram of a pseudo knot. The polynomial
\[
P_K(A, V)\ =\ (-A^{-3})^{w(K)}\, \langle K \rangle,
\]
where $w(K):=\underset{c\in C(K)}{\sum}\, sgn(c)$, $C(k)$ the set of classical crossings of $K$ and $\langle K \rangle$ the pseudo bracket polynomial of $K$, is an invariant of pseudo knots in $S^3$.
\end{theorem}

\begin{remark}
For the case of singular links in $S^3$ the reader is referred to \cite{CCC}, where three different approaches to the bracket polynomial for singular links in $S^3$ are presented.
\end{remark}

\subsection{The pseudo bracket polynomial for pseudo knots in ST}

We now generalize the pseudo bracket polynomial for pseudo links in ST. Note that we now view ST as a punctured disc (for an illustration of the looping generator $t$ in this set up see Figure~\ref{pstt}). 

\begin{figure}[ht]
\begin{center}
\includegraphics[width=3.8in]{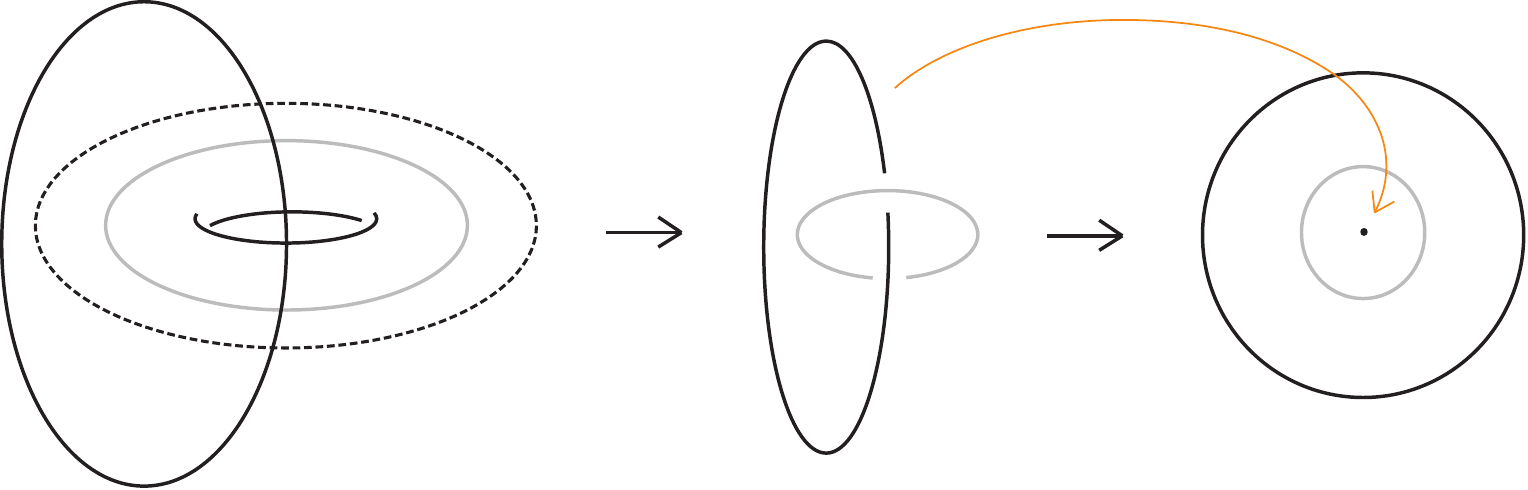}
\end{center}
\caption{The looping generator $t$ in the new set up.}
\label{pstt}
\end{figure}

We have the following:

\begin{definition}\label{pkaufbst}
Let $L$ be an oriented pseudo link in ST. The {\it pseudo bracket polynomial} of $L$ is defined by means of the relations in Definition~\ref{pkaufb} together with the following relations:
\begin{figure}[ht]
\begin{center}
\includegraphics[width=1.8in]{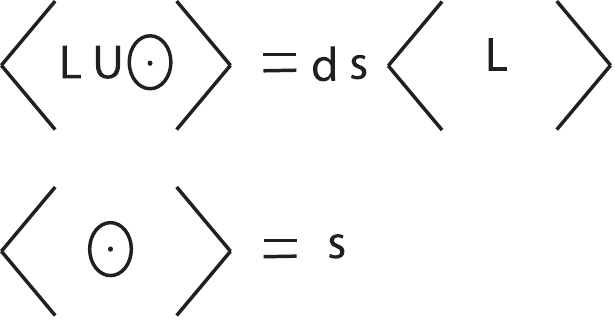}
\end{center}
%\caption{}
\label{pkbst}
\end{figure}
\end{definition}

An immediate result of Theorem~1 \cite{HD} is the following:

\begin{proposition}
The pseudo bracket polynomial for pseudo knots in ST is invariant under Reidemeister moves 2 and 3 and all pseudo moves.
\end{proposition}

As in the case of pseudo links in $S^3$ we obtain the following result:

\begin{corollary}
Let $K$ be a pseudo diagram of a pseudo knot in ST. The polynomial
\[
P_K(A, V, s)\ =\ (-A^{-3})^{w(K)}\, \langle K \rangle,
\]
where $w(K)$ is the writhe of the pseudo knots and $\langle K \rangle$ the pseudo bracket polynomial of $K$, is an invariant of pseudo knots in ST.
\end{corollary}

\begin{remark}
\begin{itemize}
\item[i.] For the case of the bracket polynomial for singular knots in ST the reader is referred to \cite{BRA}.
 
\item[ii.] For a generalization of Corollary~2 to a more general setting with an arbitrary handlebody and not just the solid torus, the reader is referred to \cite{D6}.
 
\item[iii.] It is well known that the Temperley-Lieb algebras are related to the Kauffman bracket polynomial (see \cite{LK1}). In \cite{D2}, and with the use of the Tempereley–Lieb algebra of type B, an alternative basis for the {\it Kauffman bracket skein module} of the solid torus is presented. It would be interesting to construct Kauffman-type invariants for pseudo links in ST using a generalized pseudo Temperley-Lieb algebra, and then extend this invariant to and invariant for pseudo links in the lens spaces $L(p,q)$. Note that one may also use the diagrammatic braid approach for computing Kauffman bracket skein modules (see for example \cite{D3} and \cite{D7}).
\end{itemize}
\end{remark}

\section{Conclusions}

In this paper we introduce and study pseudo links in ST via the mixed pseudo braid monoid. In particular, we present the appropriate topological set up toward the construction of HOMFLYPT-type invariants for pseudo links in $S^3$ and in ST. We believe that pseudo knots can be used in various aspects of molecular biology, and in particular, they could serve as a model to biological objects related to DNA. We also introduce and study singular links in ST and we relate the singular knot theory of ST to that of pseudo knot theory of ST. In a sequel paper we shall adopt these results in order to construct HOMFLYPT-type invariants of pseudo links in $S^3$ and ST (cor. for singular links in ST) via the generalized pseudo Hecke algebra of type B (cor. the generalized singular Hecke algebras of type B). Finally, it is worth mentioning that in \cite{D6} we study pseudo links in handlebodies and in \cite{D8} we work toward formulating the analogues of the Alexander and the Markov theorems for pseudo links in various 3-manifolds.

%%% References

\EditInfo{June 28, 2021}{February 08, 2022}{Ivan Kaygorodov}

\end{paper}